\documentclass{amsart}

\usepackage{pict2e}
\usepackage[cmtip,all]{xy}

\newtheorem{theorem}{Theorem}[section]
\newtheorem{lemma}[theorem]{Lemma}
\newtheorem{prop}[theorem]{Proposition}

\theoremstyle{definition}
\newtheorem{assumption}[theorem]{Assumption}
\newtheorem{definition}[theorem]{Definition}

\theoremstyle{remark}
\newtheorem{remark}[theorem]{Remark}

\numberwithin{equation}{section}

\newcommand{\Z}{\mathbb{Z}}
\newcommand{\Q}{\mathbb{Q}}
\newcommand{\C}{\mathbb{C}}
\newcommand{\Proj}{\mathbb{P}}
\newcommand{\calO}{\mathcal{O}}
\newcommand{\calL}{\mathcal{L}}

\newcommand{\calR}{\mathcal{R}}

\newcounter{list2}

\begin{document}
\title{Degenerations of K3 Surfaces of Degree Two}
\author{Alan Thompson}
\address{Mathematical Institute, University of Oxford, Oxford, United Kingdom, OX1 3LB.}
\curraddr{Department of Mathematical and Statistical Sciences, University of Alberta, Edmonton, Alberta, Canada, T6G 2G1.}
\email{amthomps@ualberta.ca}


\date{January 19th, 2011}

\begin{abstract} We consider a semistable degeneration of K3 surfaces, equipped with an effective divisor that defines a polarisation of degree two on a general fibre. We show that the map to the relative log canonical model of the degeneration maps every fibre to either a sextic hypersurface in $\mathbb{P}(1,1,1,3)$ or a complete intersection of degree $(2,6)$ in $\mathbb{P}(1,1,1,2,3)$. Furthermore, we find an explicit description of the hypersurfaces and complete intersections that can arise, thereby giving a full classification of the possible singular fibres. \end{abstract}

\maketitle

\section{Introduction}

The aim of this paper is to study the explicit form of the fibres arising in the relative log canonical models of semistable degenerations of K3 surfaces of degree two. This work has two main applications: compactification of the moduli of K3 surfaces of degree two and construction of explicit models for threefolds.

Degenerations of K3 surfaces were widely studied in the 1980's, with a view towards compactifying the moduli space of K3 surfaces. The first steps in this direction were taken by Kulikov \cite{kul1} \cite{kul2} and Persson-Pinkham \cite{dstcb}, who showed that a semistable degeneration of K3 surfaces may be brought into a normal form, the \emph{Kulikov model}. The central fibres in Kulikov models were classified by Persson \cite{odas}, Kulikov \cite{kul1} and Friedman-Morrison \cite{bgd1} (see Theorem \ref{kulclass}). Shepherd-Barron \cite{bgd4} then extended this work to the case of polarised degenerations, thereby showing how to study degenerations in the projective setting. The culmination of this work was a paper by Koll\'{a}r and Shepherd-Barron \cite{tadoss}, which used the fledgling minimal model programme to study the central fibres of semistable degenerations of smooth projective surfaces. There they introduce the concept of \emph{semi log canonical singularities} and show \cite[Theorem 5.1]{tadoss} that they occur as the singularities of the central fibre in the canonical model of a degeneration of surfaces of general type.  

The aim of this work is to use more recent techniques in the minimal model programme to revisit this problem in a special case: that where the K3 surfaces in question admit a polarisation of degree two. Our main result is Theorem \ref{2polclass}, which provides an explicit classification of the fibres that can occur in the relative log canonical model of such a degeneration. This classification is summarised in Tables \ref{eqtable1} and \ref{eqtable2}. We note (see Remark \ref{singrem}) that the surfaces that occur in this classification have exactly the semi log canonical singularities studied by Koll\'{a}r and Shepherd-Barron.

We remark that the problem of compactifying the moduli space of K3 surfaces of degree two has been studied before. Firstly, Shah \cite{cmsk3sd2} considered a compactification given by taking the geometric invariant theory quotient of the space of semistable (in the sense of Mumford \cite{git}) sextics in $\Proj^2$ by the action of $PGL_3$, and explicitly classified the central fibres corresponding to each of the boundary components. Later, Friedman \cite{npgttk3s} explicitly studied the boundary components arising in the Baily-Borel compactification. Whilst we consider a different semistability condition to both Shah and Friedman (arising from the work cited above and the minimal model programme), our classification is closely related to each of theirs. The precise relationship between these results is discussed in Subsection \ref{2polclasssub}.

A second application for this work is to the explicit construction of threefolds. The following abbreviated form of our main result (Theorem \ref{2polclass}) provides a classification of the canonical rings of degenerate K3 surfaces of degree two:

\begin{theorem}  Let $\pi\colon X \to \Delta$ be a semistable degeneration of K3 surfaces, with $\omega_X \cong \calO_X$. Let $H$ be a divisor on $X$ that is effective, nef and flat over $\Delta$, and suppose that $H$ induces a nef and big divisor $H_t$ on $X_t$ satisfying $H_t^2 = 2$ for $t \in \Delta^*$.

Then the morphism $\phi\colon X \to X^c$ taking $X$ to the relative log canonical model of the pair $(X,H)$ maps $X_0$ to one of:
\begin{itemize}
\item \textup{(}Hyperelliptic Case\textup{)} A sextic hypersurface 
\[\{z^2 - f_6(x_i) = 0\} \subset \Proj_{(1,1,1,3)}[x_1,x_,x_3,z].\]
\item \textup{(}Unigonal Case\textup{)} A complete intersection 
\[\{z^2 - f_6(x_i,y) = f_2(x_i) = 0\} \subset \Proj_{(1,1,1,2,3)}[x_1,x_2,x_3,y,z],\]
where $f_6(0,0,0,1) \neq 0$.
\end{itemize}
\end{theorem}

This proves an analogue of a theorem of Mendes-Lopes \cite[Theorem 3.7]{ml} that classifies the canonical rings of degenerate genus two curves. Mendes-Lopes' theorem is used by Catanese and Pignatelli \cite{flgi} to explicitly construct the relative canonical models of surfaces fibred by genus two curves. In an upcoming paper \cite{emtfk3sdt} we show that their method can be adapted to provide an explicit construction for the relative log canonical models of threefolds fibred by K3 surfaces of degree two.

\textbf{Acknowledgements.} First and foremost, I would like to thank my doctoral advisor Bal\'{a}zs Szendr\H{o}i for his support and guidance throughout the writing of this paper. I would also like to thank Miles Reid for useful conversations and for his helpful comments on an earlier draft.

\section{Background Material}\label{kulsect}

A \emph{degeneration} of K3 surfaces is a proper, flat, surjective morphism 
\begin{equation*} \pi\colon X \longrightarrow \Delta := \{z \in \C : 0 \leq |z|<1\},\end{equation*}
whose general fibre $X_t := \pi^{-1}(t)$ for $t \in \Delta^* := \Delta - \{0\}$ is a nonsingular K3 surface. Note that we do not assume anything about the algebraicity of $X$, but we do make the assumption that all components of the central fibre $X_0  := \pi^{-1}(0)$ are K\"{a}hler. 

Such a degeneration is called \emph{semistable} if $X$ is nonsingular and $X_0$ is a reduced divisor with normal crossings. By the semistable reduction theorem of Knudsen, Mumford and Waterman \cite{teI}, semistability can always be arranged via a base change and a sequence of blow-ups.  

We wish to study the relative log canonical model of a semistable degeneration of K3 surfaces. In order to do this, we first need to introduce a polarisation. A semistable degeneration $\pi\colon X \to \Delta$ is said to be \emph{polarised} if there exists a line bundle $\calL$ on $X$ that induces a nef and big line bundle $\calL_t \in \mathrm{Pic}(X_t)$ for all $t \in \Delta^*$. A theorem of Shepherd-Barron \cite[Theorem 1(a)]{bgd4} shows that, after twisting the polarisation bundle by an anti-effective divisor supported on $X_0$, we may assume that it has the form $\calL = \calO_X(H)$ for some $H$ effective and flat over $\Delta$.

In light of this, we define our main objects of study:

\begin{definition} A degeneration of K3 surfaces $\pi\colon X \to \Delta$ equipped with a divisor $H$ on $X$ is called a \emph{degeneration of K3 surfaces of degree two} if $\pi\colon X \to \Delta$ is semistable, $H$ is effective and flat over $\Delta$, and $H$ induces a nef and big divisor $H_t$ on $X_t$ satisfying $H_t^2 = 2$ for all $t \in \Delta^*$. The divisor $H$ is called the \emph{polarisation divisor} on $\pi\colon X \to \Delta$. \end{definition}

Given a degeneration of K3 surfaces of degree two $\pi\colon X\to \Delta$ with polarisation divisor $H$, in this paper we aim to study the relative log canonical model of the pair $(X,H)$ (as defined in \cite[Section 3.8]{bgav}). In order to do this, we begin by transforming the pair $(X,H)$ into a special form.

Results of Kulikov \cite{kul1} \cite{kul2} and Persson-Pinkham \cite{dstcb} show that, after a birational modification (that an easy application of \cite[Corollary 3.53]{bgav} shows does not affect the form of the relative log canonical model), we may assume that $\omega_{X} \cong \calO_{X}$. A degeneration satisfying this condition is called a \emph{Kulikov model}. We choose to work with Kulikov models because of the existence of a general classification of their central fibres. This classification was originally obtained by Persson \cite{odas}, Kulikov \cite{kul1} and Friedman-Morrison \cite{bgd1}. It states: 

\begin{theorem} \label{kulclass} \textup{\cite{bgd1}} Let $\pi\colon X \to \Delta$ be a semistable degeneration of K3 surfaces satisfying $\omega_X~\cong~\calO_X$, such that all components of $X_0$ are K\"{a}hler. Then either
\begin{list}{\textup{(\Roman{list2})}}{\usecounter{list2}}
\item $X_0$ is a smooth K3 surface;
\item $X_0$ is a chain of elliptic ruled components with rational surfaces at each end, and all double curves are smooth elliptic curves;
\item $X_0$ consists of rational surfaces meeting along rational curves which form cycles in each component. If $\Gamma$ is the dual graph of $X_0$, then $\left| \Gamma \right|$, the topological support of $\Gamma$, is homeomorphic to the sphere $S^2$.
\end{list} \end{theorem}

A Kulikov model of a degeneration of K3 surfaces will be referred to as a \emph{degeneration of Type I, II or III}, depending upon which case of the theorem it satisfies. These three cases may be distinguished by the action of the monodromy $T$ on $H^2(X_t,\Z)$ for a general fibre $X_t$. The logarithm of $T$, usually denoted by $N$, is nilpotent of index at most $3$. The possible values $(1,2,3)$ of this index correspond exactly to Types I, II and III in the theorem. 

Note that Kulikov models are not unique. Examples of this are provided by the \emph{elementary modifications of Types 0, I and II}. These are birational maps between Kulikov models that are isomorphisms outside of a codimension two subset of the central fibre $X_0$; precise definitions may be found in \cite{bgd1}. Shepherd-Barron shows \cite[Theorem 1(b)]{bgd4} that, after performing a series of elementary modifications, we may assume that the polarisation $H$ on our degeneration of K3 surfaces of degree two is nef. Furthermore, an easy application of \cite[Corollary 3.53]{bgav} shows that this process does not affect the form of the relative log canonical model.

In light of this, for the remainder of this paper we may assume that we are in the following situation: 

\begin{assumption} \label{mainass} $\pi\colon X \to \Delta$ is a degeneration of K3 surfaces of degree two with polarisation divisor $H$, such that the canonical bundle $\omega_X \cong \calO_X$ and $H$ is nef.  \end{assumption}

Under this assumption, we aim to produce an explicit classification of the central fibres that can occur in the relative log canonical model of the pair $(X,H)$.

In general, these fibres will be rather singular. We conclude this section with a brief digression to describe these singularities in more detail, before returning to classify the fibres in the next section.

More precisely, all of the singularities that will appear in the central fibre of the relative log canonical model will be \emph{semi log canonical surface singularities}. These are classified in \cite[Theorem 4.21]{tadoss}, which states that a Gorenstein surface singularity is semi log canonical if and only if it is locally analytically isomorphic to one of
\begin{itemize}
\item a smooth point;
\item a rational double point (RDP) $0 \in \{z^2 = f(x,y)\} \subset \C^3$, where the branch curve $\{f(x,y) = 0\}\subset \C^2$ has an A-D-E singularity at $0 \in \C^2$;
\item a double normal crossing point $0 \in \{xy=0\} \subset \C^3$;
\item a pinch point $0 \in \{x^2 = zy^2\} \subset \C^3$;
\item a simple elliptic singularity;
\item a cusp;
\item a degenerate cusp.
\end{itemize}

The remainder of this section will be devoted to a more detailed examination of the simple elliptic singularities, cusps and degenerate cusps. These singularities can be classified by the form of their minimal semi-resolutions (see \cite[Section 4]{tadoss} for the definition of a minimal semi-resolution).

\begin{definition}\cite[4.20]{tadoss} \label{ellsingdefn} A Gorenstein surface singularity is called:
\begin{itemize}
\item Simple elliptic if it is normal and the exceptional divisor of the minimal resolution is a smooth elliptic curve.
\item A cusp if it is normal and the exceptional divisor of the minimal resolution is a cycle of smooth rational curves or a rational nodal curve.
\item A degenerate cusp if it is not normal and the exceptional divisor of the minimal semi-resolution is a cycle of smooth rational curves or a rational nodal curve.
\end{itemize}\end{definition}

We will only be interested in some simple cases of these singularities, those which have embedding dimension $3$ (i.e. hypersurface singularities) and multiplicity $2$; it is an easy consequence of Theorem \ref{2polclass} that singularities of higher multiplicity and embedding dimension cannot occur. 

The simple elliptic singularities of embedding dimension $3$ have been classified by Saito \cite[Satz 1.9]{ees}. He finds three cases, distinguished by the self-intersection number of the exceptional elliptic curve in the minimal resolution. Only two of these cases, $\tilde{E}_7$ and $\tilde{E}_8$, have multiplicity $2$; these correspond to exceptional elliptic curves with self-intersection numbers $-2$ and $-1$ respectively. These singularities have local equations

\begin{align*} \tilde{E}_7\colon & \quad 0 \in \{z^2 = xy(y - x)(y - \lambda x)\} \subset \C^3, &\lambda \in \C-\{0,1\}, \\
\tilde{E}_8\colon & \quad 0 \in \{z^2 = y(y-x^2)(y - \lambda x^2)\} \subset \C^3, &\lambda \in \C-\{0,1\}.
\end{align*}
We will abuse notation and say that a plane curve $\{f(x,y) = 0\} \subset C^2$ has a singularity of type $\tilde{E}_7$ (resp. $\tilde{E}_8$) if the corresponding double cover $\{z^2 = f(x,y)\} \subset~\C^3$ has a singularity of type $\tilde{E}_7$ (resp. $\tilde{E}_8$).

Cusp singularities of embedding dimension $3$ have been studied by Arnold \cite{lnff}. They have the general form 
\begin{align*} T_{p,q,r}\colon & \quad 0 \in \{ x^p + y^q + z^r + \lambda xyz = 0\} \subset \C^3, & \frac{1}{p} + \frac{1}{q} + \frac{1}{r} < 1, &\quad \lambda \in \C - \{0\},\end{align*}
where the integers $p,q,r$ are determined by the form of the exceptional locus in the minimal resolution. Of these, the cusp singularities with multiplicity $2$ correspond to those $T_{p,q,r}$ with $p = 2$. Such singularities fall into two classes, $T_{2,3,r}$ with $r \geq 7$ and $T_{2,q,r}$ with $q \geq 4$ and $r \geq 5$. The forms of the exceptional loci appearing in the minimal resolutions corresponding to these cases have been calculated by Laufer in \cite[Section V]{omes}: in both cases the exceptional locus is a cycle of rational curves $E = \sum_i E_i$ with components satisfying $E_i^2 \leq -2$, that has $E^2 = -1$ (in the $T_{2,3,r}$ with $r \geq 7$ case) or $E^2 = -2$ (in the $T_{2,q,r}$ with $q \geq 4$ and $r \geq 5$ case).

In a similar way to the simple elliptic case, we will abuse notation and say that a plane curve $\{f(x,y) = 0\} \subset C^2$ has a singularity of type $T_{p,q,r}$ if the corresponding double cover $\{z^2 = f(x,y)\} \subset~\C^3$ has a singularity of type $T_{p,q,r}$.

Finally, degenerate cusps of embedding dimension $3$ and multiplicity $2$ have been classified by Shepherd-Barron in \cite[Lemma 1.3]{bgd2}. In the proof of this lemma he shows that there are two possibilities, with local equations $0 \in \{z^2 = x^2y^2\} \subset \C^3$ and $0 \in \{z^2 = y^2(y^{n}+x^2)\} \subset \mathbb{C}^3$ (where $n \geq 1$).

\section{Degenerations of K3 Surfaces of Degree Two}

We are now ready to begin our study of the central fibres occurring in the relative log canonical models of degenerations of K3 surfaces of degree two. By the results of Section \ref{kulsect}, we may assume that $\pi\colon X \to \Delta$ is a degeneration of K3 surfaces of degree two with polarisation divisor $H$ that satisfies Assumption \ref{mainass}.

Under these assumptions, it follows from the base point free theorem of Ancona \cite[Theorem 3.3]{vnvtnelbcs} that the relative log canonical algebra of the pair $(X,H)$, defined by
\begin{equation*}\calR(X,H) := \bigoplus_{n \geq 0} \pi_*\calO_X(nH),\end{equation*}
is finitely generated as an $\calO_{\Delta}$-algebra. With this in place, standard results of the minimal model programme (see \cite[Section 3.8]{bgav}, for instance) then show that the relative log canonical model of the pair $(X,H)$ exists and is given by
\begin{equation*}X^c := \mathbf{Proj}_{\Delta} \calR(X,H).\end{equation*}
The natural morphism from $X^c$ to $\Delta$ coming from this definition will be denoted $\pi^c$. Furthermore, another application of the base point free theorem shows there exists a natural birational morphism $\phi\colon X \to X^c$ over $\Delta$. A simple Riemann-Roch calculation, using the results of Mayer \cite{fk3s}, shows that this morphism takes a general fibre to either a sextic hypersurface in $\Proj(1,1,1,3)$ (the \emph{hyperelliptic case}) or a complete intersection of degree $(2,6)$ in $\Proj(1,1,1,2,3)$, where the degree two relation does not involve the degree two variable (the \emph{unigonal case}). In both cases, $\phi$ contracts at most finitely many curves in each fibre over $\Delta^*$, leading to at worst RDP singularities in the general fibre of $\pi^c\colon X^c \to \Delta$.

The aim of this paper is to find an explicit classification for the possible central fibres occurring in $X^c$. To do this we start from Theorem \ref{kulclass}, which gives a coarse classification of the possible central fibres $X_0$ of $\pi\colon X \to \Delta$ into Types I, II and III, distinguished by the index of nilpotency of the logarithm of the monodromy. We then explicitly calculate the images of these possibilities under the morphism $\phi$ to the relative log canonical model $X^c$. These images are classified by the following theorem, which is the main result of this paper.

\begin{theorem} \label{2polclass} Let $\pi\colon X \to \Delta$ be a degeneration of K3 surfaces of degree two with polarisation divisor $H$ satisfying Assumption \ref{mainass}. Then the morphism $\phi\colon X \to X^c$ taking $X$ to the relative log canonical model of the pair $(X,H)$ maps $X_0$ to one of:
\begin{itemize}
\item \textup{(}Hyperelliptic Case\textup{)} A sextic hypersurface 
\[\{z^2 - f_6(x_i) = 0\} \subset \Proj_{(1,1,1,3)}[x_1,x_,x_3,z].\]
\item \textup{(}Unigonal Case\textup{)} A complete intersection 
\[\{z^2 - f_6(x_i,y) = f_2(x_i) = 0\} \subset \Proj_{(1,1,1,2,3)}[x_1,x_2,x_3,y,z],\]
where $f_6(0,0,0,1) \neq 0$.
\end{itemize}

Furthermore, if $X_0$ is a fibre of Type I-III and $\phi(X_0)$ is hyperelliptic, then either:
\begin{list}{}{}
\item{Type I.} $\{f_6 = 0\} \subset \Proj^2$ has at worst A-D-E singularities.
\item{Type II.} $f_6 = g^2.h$ for reduced polynomials $g,h$, where either: $\mathrm{deg}(g) = 0$ and $\{h= 0\}$ has at least one singularity of type $\tilde{E}_7$ or $\tilde{E}_8$; or $\mathrm{deg}(g) > 0$ and the locus $\{g = 0\}$ is smooth and intersects $\{h = 0\}$ transversely. In either case $\{h=0\}$ may have A-D-E singularities, and if $\deg(g) >0$ it may further contain an $\tilde{E}_7$.
\item{Type III.} $f_6 = g^2.h$ for reduced polynomials $g,h$, where either: $\mathrm{deg}(g) = 0$ and $\{h= 0\}$ has exactly one singularity of type $T_{2,3,r}$ with $r \geq 7$ or type $T_{2,q,r}$ with $q \geq 4$ and $r \geq 5$; or $\mathrm{deg}(g) > 0$ and the locus $\{g = 0\}$ has either A-D-E singularities or intersects $\{h = 0\}$ non-transversely \textup{(}or both\textup{)}. In either case $\{h=0\}$ may have A-D-E singularities and cannot meet $\{g=0\}$ in any point with multiplicity $> 2$.
\end{list}
A full list of the possibilities for $\phi(X_0)$ in this case may be found in Table \ref{eqtable1}.

If $X_0$ is a fibre of Type I-III and $\phi(X_0)$ is unigonal, then:
\begin{list}{}{}
\item{Type I.} $f_2$ is irreducible and $\phi(X_0)$ has at worst RDP's.
\item{Type II.} Either $f_2$ is irreducible and $\phi(X_0)$ has at least one elliptic singularity of type $\tilde{E}_7$ or $\tilde{E}_8$; or it is a union of two distinct lines $f_2 = l_1l_2$ and the locus $\{f_6 = l_1 = l_2 = 0\} \subset \Proj(1,1,1,2)$ consists of three distinct points. In either case $\phi(X_0)$ may also contain RDP's, and in the second case it may also contain up to two $\tilde{E}_8$'s.
\item{Type III.} Either $f_2$ is irreducible and $\phi(X_0)$ has exactly one cusp singularity of type $T_{2,3,r}$ with $r \geq 7$ or type $T_{2,q,r}$ with $q \geq 4$ and $r \geq 5$; or it is a union of two distinct lines $f_2 = l_1l_2$ and the locus $\{f_6 = l_1 = l_2 = 0\} \subset \Proj(1,1,1,2)$ consists of two points, one of which has multiplicity two. In either case $\phi(X_0)$ may also contain some RDP's.
\end{list}
A full list of the possibilities for $\phi(X_0)$ in this case may be found in Table \ref{eqtable2}.
\end{theorem}

\begin{table}
\caption{Possibilities for $\phi(X_0) = \{z^2 - f_6(x_i) = 0\} \subset \Proj(1,1,1,3)$ hyperelliptic.}
\label{eqtable1}
\begin{tabular}{|c|c|c|l|}
\hline
Type & Name & $f_6(x_i)$ & Comments \\
\hline
I & h & Reduced & $f_6$ has at worst A-D-E's.\\
\hline
II & 0h & Reduced & $f_6$ has one $\tilde{E}_7$, one $\tilde{E_8}$ or two $\tilde{E_8}$'s.\\
& 1 & $l^2(x_i)f_4(x_i)$ & $l$ linear, $|l\cap f_4| = 4$, where $f_4$ may have an $\tilde{E}_7$.\\
& 2 & $q^2(x_i)f_2(x_i)$ & $q$ smooth quadric, $|q \cap f_2|= 4$.\\
& 3 & $f_3^2(x_i)$ & $f_3$ smooth cubic. \\
\hline
III & 0h & Reduced & $f_6$ has exactly one $T_{2,3,r}$ with $r \geq 7$ or $T_{2,q,r}$ \\ 
&&& with $q \geq 4$ and $r \geq 5$.\\
& 1 & $l^2(x_i)f_4(x_i)$ & $l$ linear, $|l\cap f_4| \leq 3$ with multiplicities $\leq 2$.\\ 
& 2 & $q^2(x_i)f_2(x_i)$ & $q$ (possibly nodal) quadric, $|q \cap f_2| \leq 4$ \\ 
&&& ($<4$ if $q$ smooth) with multiplicities $\leq 2$.\\
& 3 & $f_3^2(x_i)$ & $f_3$ cubic with nodal singularities. \\
\hline
\end{tabular}
\end{table}

\begin{table}
\caption{Possibilities for $\phi(X_0) = \{z^2 - f_6(x_i,y) = f_2(x_i) = 0\} \subset \Proj{(1,1,1,2,3)}$ unigonal.}
\label{eqtable2}
\begin{tabular}{|c|c|c|l|}
\hline
Type & Name & $f_2(x_i)$ & Comments \\
\hline
I & u & Irreducible & $\phi(X_0)$ has at worst RDP's.\\
\hline
II & 0u & Irreducible & $\phi(X_0)$ has one $\tilde{E}_7$, one $\tilde{E_8}$ or two $\tilde{E_8}$'s.\\
& 4 & $l_1(x_i)l_2(x_i)$ & $l_i$ linear, $|l_1 \cap l_2 \cap f_6| = 3$, where $\phi(X_0)$ may \\
&&& have one or two $\tilde{E}_8$'s.\\
\hline
III & 0u & Irreducible & $\phi(X_0)$ has exactly one $T_{2,3,r}$ with $r \geq 7$ or $T_{2,q,r}$\\
&&& with $q \geq 4$ and $r \geq 5$.\\ 
& 4 & $l_1(x_i)l_2(x_i)$ & $l_i$ linear, $|l_1 \cap l_2 \cap f_6| = 2$, where the curve \\
&&& $\{f_6 = l_i = 0\}$ may be non-reduced for one or \\ 
&&& both choices of $i \in \{1,2\}$.\\
\hline
\end{tabular}
\end{table}

\begin{remark} \label{singrem} We remark that $\phi(X_0)$ has semi log canonical singularities in all cases of this theorem. Furthermore, all of the semi log canonical singularities of embedding dimension $3$ and multiplicity $2$ listed in Section \ref{kulsect} occur and the classification of these singularities given there can be deduced from the statement of the theorem.\end{remark}

\subsection{Relation to Other Results}\label{2polclasssub} 

We use this subsection to compare the cases in this theorem to the description of the boundary components in two known compactifications of the moduli space of K3 surfaces of degree two. 

Firstly, in \cite[Section 5]{npgttk3s}, Friedman obtains arithmetic possibilities for the the four Type II boundary components appearing in the Baily-Borel compactification of the moduli of K3 surfaces of degree two. This compactification is quite simple to describe, with boundary components determined by the form of the monodromy weight filtration on $H^2(X_t,\Q)$. Due to this simplicity, we expect to have a good match between the cases in Friedman's classification and those in Tables \ref{eqtable1} and \ref{eqtable2}. The exact nature of this correspondence is described in Table \ref{friedtab}.
\begin{table}
\caption{Correspondence between cases in \cite[Section 5]{npgttk3s} and Type II fibres in Tables \ref{eqtable1} and \ref{eqtable2}.}
\label{friedtab}
\begin{tabular}{|c|c|l|}
\hline
Case in \cite[Section 5]{npgttk3s} & Case in Tables \ref{eqtable1} and \ref{eqtable2} & Comments\\
\hline
(5.2.1) & (II.3) &\\
(5.2.2) & (II.0h), (II.0u) or (II.4) & (II.0h), (II.0u) contain an $\tilde{E}_8$\\
(5.2.3) & (II.0h), (II.0u) or (II.1) & (II.0h), (II.0u) contain an $\tilde{E}_7$ \\
(5.2.4) & (II.2) &\\
\hline
\end{tabular}
\end{table}

Note here that cases (II.0h) and (II.0u) from Tables \ref{eqtable1} and \ref{eqtable2} do not appear as separate cases in Friedman's list. This is because the Baily-Borel compactification depends only upon the monodromy weight filtration and not upon the form of the polarisation on the central fibre, so Friedman is free to twist his polarisation by divisors supported on $X_0$. This allows him \cite[Theorem 2.2]{npgttk3s} to assume that $H_0.D_i > 0$ for all double curves $D_i \subset X_0$, which eliminates these two cases from consideration. However, in order to study the relative log canonical models of degenerations of K3 surfaces of degree two, in our case we wish to degenerate pairs consisting of a K3 surface of degree two along with its polarisation divisor, so we are not free to perform such twisting operations (as they would alter the polarisation). This gives rise to cases (II.0h) and (II.0u), which correspond to Type II degenerations where $H_0.D_i = 0$ for all double curves $D_i \subset X_0$. In the Baily-Borel compactification these cases appear in two different strata, distinguished by the types of singularities occurring in them; this correspondence is detailed in Table~\ref{friedtab}.

Secondly, in \cite{cmsk3sd2}, Shah finds a compactification for the moduli space of K3 surfaces of degree two, by taking the geometric invariant theory quotient of the space of semistable (in the sense of Mumford \cite{git}) sextics in $\Proj^2$ by the action of $PGL_3$. In \cite[Theorem 2.4]{cmsk3sd2} he classifies the boundary components of this compactification by classifying the semistable sextics corresponding to closed $PGL_3$-orbits. As was the case with Friedman's classification, there is a matching between these cases and the cases from Tables \ref{eqtable1} and \ref{eqtable2}, summarised in Table \ref{shahtab}. Note however, that this matching is not quite as neat as in Friedman's case. This occurs because of the difference between our semistability condition and Shah's, meaning that the cases in Tables \ref{eqtable1} and \ref{eqtable2} may comprise a union of (not necessarily closed) semistable $PGL_3$-orbits.
\begin{table}
\caption{Correspondence between cases in \cite[Theorem 2.4]{cmsk3sd2} and in Tables \ref{eqtable1} and \ref{eqtable2}.}
\label{shahtab}
\begin{tabular}{|c|c|l|}
\hline
Case in \cite[Theorem 2.4]{cmsk3sd2} & Case in Tables \ref{eqtable1} and \ref{eqtable2} & Comments\\
\hline
(II.1) & (II.0h) & (II.0h) contains an $\tilde{E}_8$. \\
(II.2) & (II.0h) or (II.1) & (II.0h) contains an $\tilde{E}_7$.\\
(II.3) & (II.2) & \\
(II.4) & (II.3) & \\
\hline
(III.1) & (III.0h) or (III.2) & (III.0h) contains a $T_{2,3,r}$,\\
& & with $r \geq 7$.\\
(III.2) & (III.0h), (III.1) or (III.3) & (III.0h) contains a $T_{2,q,r}$,\\
& & with $q \geq 4$ and $r \geq 5$.\\
\hline
(IV) & (II.0u), (II.4), & \\
& (III.0u) or (III.4) & \\
\hline
\end{tabular}
\end{table}

Note also that all of the unigonal fibres map to a single point in Shah's compactification, represented by a fibre in group (IV). In order to solve this problem, in \cite[Section 4]{cmsk3sd2} Shah resolves this point to separate out the unigonal cases. He obtains \cite[Theorem 4.3]{cmsk3sd2}, which classifies these fibres. The comparison with the unigonal cases in our theorem can be found in Table \ref{shahtab2}.
\begin{table}
\caption{Correspondence between cases in \cite[Theorem 4.3]{cmsk3sd2} and in Tables \ref{eqtable1} and \ref{eqtable2}.}
\label{shahtab2}
\begin{tabular}{|c|c|}
\hline
Case in \cite[Theorem 4.3]{cmsk3sd2} & Case in Tables \ref{eqtable1} and \ref{eqtable2}\\
\hline
1(ii) & (II.0u) \\
2(i) & (II.4)\\
2(ii) & (III.0u) or (III.4) \\
\hline
\end{tabular}
\end{table}

\section{Components of Degenerate Fibres} \label{compdenfib} 

In order to prove Theorem \ref{2polclass}, we need to find a way to study the effects of the morphism $\phi$ on $X_0$. To do this, we begin by showing that, under the assumptions of Theorem \ref{2polclass}, the log canonical model 
\begin{equation*}(X_0)^c := \mathrm{Proj} \bigoplus_{n\geq 0} H^0(X_0, \calO_{X_0}(nH_0))\end{equation*}
of the pair $(X_0,H_0)$, where $H_0$ denotes the divisor induced on $X_0$ by $H$, agrees with the central fibre $(X^c)_0$ of $\pi^c\colon X^c \to \Delta$. We can then use the linear systems induced by $H_0$ on the components of $X_0$ to study the natural morphism $\phi_0\colon X_0 \to (X_0)^c$. As $\phi_0$ agrees with the restriction of $\phi$ to $X_0$, this will enable us to prove the theorem.

Our first step is to show that $(X^c)_0$ and $(X_0)^c$ agree. As noted above, the log canonical model $(X_0)^c$ of the pair $(X_0, H_0)$ is defined by the global sections $H^0(X_0,\calO_{X_0}(nH_0))$ for $n>0$. On the other hand, the central fibre $(X^c)_0$ of the relative log canonical model $\pi^c\colon X^c \to \Delta$ is defined by the localised direct images $\pi_*(\calO_X(nH))_0 \otimes_{\calO_{\Delta,0}} k(0)$ for $n>0$, where $k(0)$ is the residue field at $0 \in \Delta$. The fact that these two maps agree follows immediately from the following lemma, which is a slight generalisation of \cite[Lemma 2.17]{bgd4}.

\begin{lemma} \label{fibreres} Suppose that $\pi\colon X \to \Delta$ is a degeneration of K3 surfaces of degree two with polarisation divisor $H$ that satisfies Assumption \ref{mainass}. Then the natural maps
\begin{equation*}\pi_*(\calO_X(nH))_0 \otimes_{\calO_{\Delta,0}} k(0) \longrightarrow H^0(X_0, \calO_{X_0}(nH_0))\end{equation*}
are isomorphisms for all $n>0$.
\end{lemma}
\begin{proof} By \cite[Theorem 2.1]{vnvtnelbcs}, the assumptions on $X$ and $H$ imply that the higher direct images $R^i\pi_*(\calO_X(nH)) = 0$ for all $i>0$ and all $n>0$. The result then follows easily from the Theorem on Cohomology and Base Change.\end{proof}

This result enables us to restrict our attention to the pair $(X_0,H_0)$. We use the remainder of this section to collect together some results on the interaction of $H_0$ with the components of $X_0$, which will come in useful in the proof of Theorem \ref{2polclass}.

We begin by fixing some notation. Suppose that $\pi\colon X \to \Delta$ is a degeneration of K3 surfaces of degree two with polarisation divisor $H$ that satisfies Assumption \ref{mainass}, with central fibre $X_0 = \pi^{-1}(0)$ of Type II or III. By the classification of the central fibres of Kulikov models (Theorem \ref{kulclass}), $X_0$ is a union of rational and elliptic ruled components meeting transversely along a set of double curves. Write $X_0 = V_1 \cup \cdots \cup V_r$ where the $V_i$ are the irreducible components of $X_0$ and we assume that the $V_i$ have been normalised. Let $D_{ij}$ denote the double curve $V_i \cap V_j$ and let $D_i = \bigcup_j D_{ij}$ denote the double locus on $V_i$. Let $H_i$ denote the effective (or zero) divisor obtained by restricting $H_0$ to the component $V_i$.

To study the behaviour of the polarisation on $V_i$, we follow Shepherd-Barron \cite{bgd4} and start by separating the $V_i$ into three sets according to the properties of $H_i$. We will call a curve $C \subset V_i$ a \emph{$0$-curve} if $H_i . C = 0$. Then $V_i$ will be called a \emph{$0$-surface} if it contains only finitely many $0$-curves; a \emph{$2$-surface} if $H_i$ is numerically trivial; and a \emph{$1$-surface} if it contains a pencil of $0$-curves but is not a $2$-surface. Note that these classes are mutually exclusive and that, by \cite[Proposition 2.3]{bgd4}, every component of $X_0$ is either a $0$-, $1$- or $2$-surface.

This classification will be useful because, as we shall see later, the map $\phi$ to the relative log canonical model of the pair $(X,H)$ defines a birational morphism on each $0$-surface, contracts each $1$-surface to a curve and contracts each $2$-surface to a point. This observation will allow us to calculate the possible images of $X_0$ under $\phi$ by studying the possible configurations of $0$-, $1$- and $2$-surfaces that can occur in it.

The following result about $0$-, $1$- and $2$-surfaces is an easy consequence of \cite[Proposition 2.3]{bgd4}:

\begin{prop} \label{012facts} Suppose that $V_i$ and $H_i$ are defined as above. Then
\begin{list}{\textup{(\roman{list2})}}{\usecounter{list2}}
\item If $V_i$ is a $1$-surface, then the pencil of $0$-curves on $V_i$ forms a ruling.
\item If $V_i$ contains a pencil of $0$-curves and another $0$-curve that does not lie in this pencil, then $V_i$ is a $2$-surface.
\item If $V_i$ contains an effective divisor $E$ that satisfies $E^2 > 0$ and $H_i.E =0$, then $V_i$ is a $2$-surface.
\end{list} \end{prop}
\begin{proof} (i) and (ii) are immediate from \cite[Proposition 2.3]{bgd4}. (iii) follows easily from the Hodge Index Theorem.\end{proof}

In addition to this, we have the following information about $0$- and $1$-surfaces:

\begin{lemma} \label{Hpos} $V_i$ is a $0$-surface if and only if $H_i^2 >0$. Furthermore, $V_i$ is a $1$-surface if and only if $H_i^2 = 0$ and $H_i.D_i > 0$. \end{lemma}
\begin{proof} First assume that $V_i$ is a $0$-surface. Then $H_i^2 >0$ by \cite[Lemma 2.8]{bgd4}.

Next, assume that $H_i^2 > 0$. Then $H_i$ cannot be numerically trivial, so $V_i$ is not a $2$-surface. So suppose that $V_i$ is a $1$-surface. Then, by Proposition \ref{012facts}, $V_i$ is ruled and the pencil of $0$-curves form a ruling. Let $F$ be any such $0$-curve. Then $F^2 = 0$ and, since $H_i^2 > 0$, the Hodge Index Theorem implies that $F$ is numerically equivalent to 0. But $F$ is a fibre of a ruling, so this cannot occur. Thus, $V_i$ is not a $1$- or $2$-surface, so it must be a $0$-surface. 

Now suppose that $H_i^2 = 0$ and $H_i.D_i > 0$. The argument above then shows that $V_i$ is not a $0$-surface and, as $H_i$ is not numerically trivial, $V_i$ cannot be a $2$-surface either. So $V_i$ is a $1$-surface.

Finally, assume that $V_i$ is a $1$-surface. Then the argument above shows that $H_i^2 = 0$. Furthermore, by Proposition \ref{012facts}, the $0$-curves form a ruling on $V$. Let $F$ denote an irreducible fibre of this ruling. Then, by adjunction, $F.D_i = 2$. So $D_i$ must contain an irreducible component $D_{ij}$ that is not contained in a fibre of the ruling on $V_i$ and thus, by Proposition \ref{012facts} again, must satisfy $H_i.D_{ij} >0$. So, as $H_i$ is nef, $H_i.D_i > 0$ also. \end{proof}

For the rest of this section, we must separate the cases where the components under consideration are rational or elliptic ruled. This distinction will enable us to get much more information about the components themselves and the polarisation divisors on them.

So suppose first that $V_i$ is a rational component with locus of double curves $D_i$. Then $(V_i, D_i)$ is an example of an \emph{anticanonical pair}: a pair $(V,D)$ consisting of a rational surface $V$ and a reduced section $D \in \left| -K_{V} \right|$ (which is necessarily connected). Anticanonical pairs have been studied extensively by Friedman in \cite{bgd4a}; all of the results about them that we will need may be found in his paper and so will not be reproduced here.

Thus, we are left with the case where $V_i$ is elliptic ruled. In the following subsection we will perform a brief study of such surfaces, culminating in an analogue of a result of Friedman \cite[Theorem 10]{bgd4a} describing the properties of certain linear systems on them. 

\subsection{Elliptic Ruled Components}\label{ellcomp}

By the classification of Kulikov models (Theorem \ref{kulclass}), elliptic ruled components can only appear in a degeneration of Type II. Furthermore, each elliptic ruled component in a Type II degeneration contains precisely two smooth elliptic double curves that form sections for the ruling. With this in mind, define:

\begin{definition} An \emph{anticanonical triple} $(V,D',D'')$ is a (not necessarily minimal) elliptic ruled surface $V$, along with two disjoint smooth elliptic curves $D'$ and $D''$ that form sections for the ruling and satisfy $K_V \sim -D'-D''$.\end{definition}

In this subsection we aim to emulate some of Friedman's \cite{bgd4a} results for anticanonical pairs in the new setting of anticanonical triples. We first add a polarisation divisor.

\begin{definition} Let $(V,D',D'')$ be an anticanonical triple and suppose that $H$ is an effective and nef divisor on $V$. Then $(V,D',D'')$ is called \emph{$H$-minimal} if there do not exist any rational $(-1)$-curves $C$ on $V$ with $H.C = 0$. \end{definition}

With this in place, we begin our study with a result describing the interaction between the polarisation divisor and the double curves.

\begin{lemma}\label{HD'0} Let $(V,D',D'')$ be an anticanonical triple and let $H$ be an effective and nef divisor on $V$ with $0< H^2 \leq 2$. Suppose that $(V,D',D'')$ is $H$-minimal. Then:
\begin{list}{\textup{(\roman{list2})}}{\usecounter{list2}}
\item If $(D')^2 = (D'')^2 = 0$, then $V$ is minimally ruled and $H.D' = H.D'' > 0$.
\item If $(D')^2$ and $(D'')^2$ are not both zero, then at least one of $D'$ and $D''$ ($D'$ say) satisfies $-H^2 \leq (D')^2 < 0$ and $H.D' =0$.
\end{list} \end{lemma}
\begin{proof} Consider (i) first. By \cite[Proposition III.18]{cas}, if $V$ is minimally ruled then $(D')^2 + (D'')^2 = 0$. Each time $V$ is blown up this number decreases by $1$, so $(D')^2 + (D'')^2 \leq 0$ with equality if and only if $V$ is minimally ruled. So if $(D')^2$ and $(D'')^2$ are both zero then $V$ is minimally ruled. By \cite[Proposition III.18]{cas} again, we see that in this case $D' \sim D''$, so $H.D' = H.D''$. Finally, this intersection number is strictly positive by the Hodge index theorem.

Now we prove (ii). As $(D')^2 + (D'')^2 \leq 0$, if $(D')^2$ and $(D'')^2$ are not both zero then at least one of $(D')^2$ and $(D'')^2$ must be strictly negative. We split into two cases: that where one of $(D')^2$ and $(D'')^2$ is strictly negative and the other is not, and that where both $(D')^2$ and $(D'')^2$ are strictly negative.

In the first case, without loss of generality we can assume that $(D')^2<0$ and $(D'')^2~\geq~0$. In this case, analysis of $H^2(V,\Z)$ using \cite[Proposition III.18]{cas} shows that $D'' - D'$ is numerically equivalent to an effective sum of components of fibres of the ruling on $V$. So, as $H$ is nef, $H.(D''-D') \geq 0$ and $H.D'' \geq H.D'$.

In the second case, as both $(D')^2$ and $(D'')^2$ are strictly negative, without loss of generality we may choose $D'$ so that $H.D' \leq H.D''$.

In both cases, we see that we can choose $D'$ so that $(D')^2<0$ and $H.D' \leq H.D''$. Now, noting that $H^2(V, \calO_V(H-D')) \cong H^0(V, \calO_V(-H-D''))$ by Serre duality and that this second group vanishes as $(H+D'')$ is effective, by Riemann-Roch
\begin{equation*} h^0(V,\calO_V(H-D'))\geq \frac{1}{2}(H^2-H.D'+H.D'') > 0.\end{equation*}
So $H-D' \sim E$ for some effective divisor $E$.

We next show that $E$ is nef. Suppose for a contradiction that there exists $C$ irreducible such that $E.C<0$. Note that $C^2<0$, as $E$ is effective. If $C = D'$ then $E.C = H.D' - (D')^2 >0$ as $(D')^2 <0$ and $H$ is nef, contradicting $E.C<0$. So $C \neq D'$ and $C.D'\geq 0$. If $C.D' = 0$, then $E.C = H.C \geq 0$, as $H$ is nef, again contradicting $E.C < 0$. So $C.D'>0$. In this case, the genus formula gives $C^2 = -1$ and $C.D'=1$. But then $0>E.C = H.C - 1 $, so $H.C = 0$. Thus, $C$ must be a rational $(-1)$-curve with $H.C=0$, which cannot exist by the $H$-minimality assumption. Thus, no such $C$ may exist and $E$ is nef.

Therefore $0 \leq E^2= H^2 -2H.D' + (D')^2$. But $H^2 \leq 2$ and $(D')^2 < 0$, so we must have $H.D' = 0$ and $-H^2 \leq (D')^2$.
\end{proof}

Next we have a result, analogous to \cite[Lemma 5]{bgd4a}, that will allow us to calculate the dimension of the linear system induced by the polarisation in certain cases.

\begin{lemma} \label{RRtriple} Let $(V,D',D'')$ be an anticanonical triple and let $H$ be an effective divisor on $V$. Then
\begin{list}{\textup{(\roman{list2})}}{\usecounter{list2}}
\item $h^0(V,\calO_V(H)) - h^1(V,\calO_V(H)) = \frac{1}{2}(H^2 + H.D' + H.D'').$ 
\end{list}
Furthermore, suppose that $H$ is nef with $H^2>0$. If $H.D' = 0$ and $D'$ is not fixed in $|H|$, then:
\begin{list}{\textup{(\roman{list2})}}{\usecounter{list2}} \addtocounter{list2}{1}
\item If $H.D'' >0$ then $h^1(V,\calO_V(H)) = 1$.
\item If $H.D'' = 0$ and $D''$ is not fixed in $|H|$, then $h^1(V,\calO_V(H)) = 2$.
\end{list} \end{lemma}
\begin{proof} By Serre duality, $H^2(V,\calO_V(H)) \cong H^0(V,\calO_V(-H-D'-D''))$, which vanishes as $(H+D'+D'')$ is effective. Given this, (i) follows immediately from the Riemann-Roch theorem. Statements (ii) and (iii) follow easily from the exact sequence
\begin{equation*} 0 \longrightarrow \calO_V(H-D'-D'') \longrightarrow \calO_V(H) \longrightarrow \calO_{D'+D''}(H) \longrightarrow 0.\end{equation*}\end{proof}

Finally, we are now in a position to prove an analogue of Friedman's main result \cite[Theorem 10]{bgd4a} describing the behaviour of linear systems on anticanonical triples.

\begin{theorem} \label{actriple} Let $(V,D',D'')$ be an anticanonical triple and let $H$ be an effective and nef divisor on $V$ with $H^2>0$. Suppose that $(V,D',D'')$ is $H$-minimal, neither $D'$ nor $D''$ is fixed in $|H|$ and $H.D'=0$. Denote the mobile part of the linear system $|H|$ by $|H|_m$ and the fixed part by $|H|_f$, so that $|H| = |H|_m + |H|_f$. Then
\begin{list}{\textup{(\roman{list2})}}{\usecounter{list2}}
\item If $|H|_f = 0$, then $|H|$ has base points if and only if $H.D'' = 1$, in which case $p = H.D''$ is the only base point and the general member of $|H|$ is smooth.
\item If $|H|_f \neq 0$, then $|H|_f$ is a smooth rational curve $F$ with $F^2 = -1$ or $-2$ and $|H|_m = kE$ for $k\geq 2$ and $E$ smooth elliptic with $E^2 = E.D'' = 0$ and $F.E = 1$.
\item $|nH|$ has no fixed components or base locus for $n\geq 2$. If $n \geq 3$, the corresponding morphism $\phi_{|nH|}$ is birational onto its image.
\end{list}
\end{theorem} 

\begin{remark} \label{ellrem} We remark that this theorem depends upon the assumption that $H.D' = 0$ which, by Lemma \ref{HD'0}, is always true when $H^2 \leq 2$. For larger values of $H^2$, however, this assumption will not necessarily hold. That said, we expect that a similar, albeit substantially more complicated, result should hold without the $H.D'=0$ assumption, which would allow generalisation to larger values of $H^2$. \end{remark}

\begin{proof}[Proof of Theorem \ref{actriple}] (Based upon the proof of \cite[Theorem 10]{bgd4a}). We prove (i) first. Suppose that $H.D'' = 0$. Then if $Z$ is any reduced member of $|H|$, adjunction gives $H|_Z \cong \omega_Z$. Thus, by
\begin{equation*} 0 \longrightarrow \calO_V \longrightarrow \calO_V(H) \longrightarrow \omega_Z \longrightarrow 0.\end{equation*}
and Lemma \ref{RRtriple}(iii) we see that $H^0(V,\calO_V(H)) \to H^0(Z,\omega_Z)$ is surjective. So, by \cite[Lemma 2]{bgd4a}, $|H|$ is base point free.

Next suppose that $H.D''\geq 2$. The exact sequence
\begin{equation*} 0 \longrightarrow \calO_V(H-D'-D'') \longrightarrow \calO_V(H-D'') \longrightarrow \calO_{D'}(H-D'') \longrightarrow 0.\end{equation*}
gives $h^1(V,\calO_V(H-D''))=1$, so the exact sequence
\begin{equation*} 0 \longrightarrow \calO_V(H-D'') \longrightarrow \calO_V(H) \longrightarrow \calO_{D''}(H) \longrightarrow 0.\end{equation*}
and Lemma \ref{RRtriple}(ii) show that $H^0(V,\calO_V(H)) \to H^0(D'',\calO_{D''}(H))$ is surjective. But $h^0(D'',\calO_{D''}(H)) = \deg_{D''}(H) \geq 2$, so $|H|$ cannot have base points on $D''$. Given this, the remainder of the proof of (i) proceeds exactly as the proof of \cite[Theorem 10.2]{bgd4a}.

Before embarking upon the proof of (ii), we prove a short lemma:

\begin{lemma} \label{Hflem} Under the assumptions of Theorem \ref{actriple}, if $|H|_m.|H|_f = 0$ then $|H|_f=0$. \end{lemma}
\begin{proof} Assume that $|H|_m.|H|_f = 0$. If $F$ is any component of $|H|_f$, then we have $|H|_f.F = H.F \geq 0$ so, as $|H|_f$ is effective, $|H|_f$ must be nef. Thus, if we also have $|H|_f^2 > 0$, Lemma \ref{RRtriple} applied to $|H|_f$ shows that $h^0(V,\calO_V(|H|_f)) >1$, contradicting $|H|_f$ being fixed. 

Therefore we may assume $|H|_f^2 = 0$, which implies that $H.|H|_f =0$. So, as $H$ is nef, any component $F$ of $|H|_f$ must satisfy $H.F=0$ and, by the Hodge Index Theorem, $F^2<0$. Furthermore, by $H$-minimality, $F^2=-1$ cannot occur. So, by the genus formula, $F$ must be a rational $(-2)$-curve. Therefore $|H|_f$ is an effective sum of rational $(-2)$-curves, which occur in chains supported on fibres of the ruling on $V$. However, analysis of such configurations shows that $|H|_f^2 = 0$ is impossible unless $|H|_f~=~0$. \end{proof}

Now we prove (ii). Suppose first that $|H|_m^2>0$. By applying Lemma \ref{RRtriple}(i) to $H$ and $|H|_m$ and subtracting (noting that $h^0(V,\calO_V(H)) =h^0(V,\calO_V(|H|_m))$), we obtain
\begin{equation}  h^1(V,\calO_V(|H|_m)) - h^1(V,\calO_V(H)) = \frac{1}{2}(|H|_m.|H|_f + H.|H|_f + |H|_f.D'').\label{H1eq}\end{equation}
Since $H$ and $|H|_m$ are nef and $D''$ is not a component of $|H|_f$, all terms on the right hand side are $\geq 0$. Furthermore, if $h^1(V,\calO_V(H)) = h^1(V,\calO_V(|H|_m)$  we must have $|H|_m.|H|_f =0$, so in this case $|H|_f=0$ by Lemma \ref{Hflem}. By Lemma \ref{RRtriple}, the only remaining case is $h^1(V,\calO_V(H)) = 1$ and  $h^1(V,\calO_V(|H|_m) = 2$, occurring when $H.D''>0$ and $|H|_m.D'' =0$. In this case we must have $|H|_f.D'' > 0$, and we may assume $|H|_m.|H|_f >0$ as otherwise we would have $|H|_f=0$ by Lemma \ref{Hflem}. Using this, Equation \eqref{H1eq} gives $|H|_f.D'' = |H|_m.|H|_f =1$ and $H.|H|_f = 0$ and, by the Hodge index theorem, all components of $|H|_f$ must have strictly negative self-intersection. Under these conditions, the argument used to prove \cite[Theorem 10.1.2]{bgd4a} holds to show that $|H|_f$ must contain a rational $(-1)$-curve $C$ with $H.C=0$, contradicting $H$-minimality.

Thus we are left with the case where $|H|_m^2=0$ and $|H|_m.|H|_f>0$. In this case $|H|_m$ is base point free (by part (i)) and has no fixed components so, by Bertini's theorem, either $|H|_m$ contains a smooth irreducible member or $|H|_m$ is composite with a pencil.

Let $Z$ be an irreducible component of $|H|_m$. By the genus formula $Z$ is either rational with $Z.D'' = 2$ or elliptic with $Z.D'' =0$, and the sequence
\begin{equation*} 0 \longrightarrow \calO_V \longrightarrow \calO_V(Z) \longrightarrow \calO_Z(Z) \longrightarrow 0.\end{equation*}
shows that $h^0(V,\calO_V(Z)) = 2$.  

If $Z$ is rational, the exact sequence above gives $h^1(V,\calO_V(Z)) = 1$ and induction on the number of components of $|H|_m$ gives $h^1(V,\calO_V(|H|_m)) = 1$. But then equation \eqref{H1eq} gives that $|H|_m.|H|_f = 0$, which is a contradiction.

So $Z$ is elliptic with $Z.D'' =0$ and if $|H|_m$ is composite with a pencil, this pencil is necessarily rational. So $|H|_m = kZ$ for some $k\geq 1$. The exact sequence above gives $h^1(V,\calO_V(Z)) = 2$. 
 
Now, as $|H|_m.|H|_f>0$, there exists an irreducible component $F$ of $|H|_f$ with $Z.F >0$. By Lemma \ref{RRtriple} applied to $F$, we obtain $F^2 \leq 0$. If $F^2 = 0$ then $(Z+F)$ is nef with $(Z+F)^2 > 0$. So in a manner analogous to equation \eqref{H1eq} we obtain
\begin{equation*}  h^1(V,\calO_V(Z)) - h^1(V,\calO_V(Z+F)) = \frac{1}{2}(Z.F + (Z+F).F + F.D'').\end{equation*}
Analysing this using Lemma \ref{RRtriple}, we see that $(Z+F).F = 0$. But then as $F^2=0$ we must have $Z.F=0$, a contradiction. So $F^2<0$ and by the genus formula $F$ is smooth and rational with $F^2= -1$ or $-2$.

Suppose that $Z.F = a >1$. Assume that $F^2 = -2$. Then $(Z+F)$ is nef and has $(Z+F)^2>0$, so applying Lemma \ref{RRtriple}(iii) (noting that $F.D'' = 0$) to $(Z+F)$ we obtain $h^1(V,\calO_V(Z+F)) = 2$. Using this, the sequence
\begin{equation*} 0 \longrightarrow \calO_V(Z) \longrightarrow \calO_V(Z+F) \longrightarrow \calO_F(Z+F) \longrightarrow 0.\end{equation*}
shows that $H^0(V,\calO_V(Z+F)) \to H^0(F,\calO_F(Z+F))$ is surjective. So $F$ is not fixed in $(Z+F)$ and hence cannot be in $|H|_f$, giving a contradiction. A similar argument in the case $F^2 = -1$ shows that $H^0(V,\calO_V(Z+F)) \to H^0(F,\calO_F(Z+F))$ has image of codimension $1$, giving the same contradiction.

Therefore $Z.F = 1$. If $k=1$, the argument used to prove \cite[Theorem 10.1.2]{bgd4a} holds here to show that $|H|_f$ contains a rational $(-1)$-curve $C$ with $H.C = 0$, contradicting $H$-minimality. So we must have $k\geq 2$.

To complete the proof of (ii), suppose first that $H.D'' = 0$. Let $H' = kZ+F$, where $F$ as above satisfies $F^2 = -2$ and $Z.F = 1$ ($F^2=-1$ cannot occur as $H.D'' =0$). Define $H'' = H-H'$. Then $H'$ is nef and has $(H')^2 > 0$. By Lemma \ref{RRtriple}(iii), we see that $h^1(V,\calO_V(H)) = h^1(V,\calO_V(H')) =2$, and a calculation analogous to that used to obtain equation \eqref{H1eq} gives $0 = H'.H'' + H.H''$. But then $H'.H'' = 0$ and the argument used to prove Lemma \ref{Hflem} holds here to show that $H'' = 0$. Thus, $H = kZ+F$, proving one case of (ii).

Finally, suppose that $H.D'' >0$. In this case, if $F^2 = -2$ the argument used to prove \cite[Theorem 10.1.2]{bgd4a} holds here to show that $|H|_f$ contains a rational $(-1)$-curve $C$ with $H.C = 0$, contradicting $H$-minimality. If $F^2 = -1$ then $F.D''=1$ by the genus formula, and an argument analogous to the one above shows that $H = kZ+F$. This completes the proof of (ii).

It just remains to prove (iii). If $H.D''>0$, the proof of (iii) proceeds exactly as the proof of \cite[Theorem 10.4]{bgd4a}. The remaining case is $H.D''=0$. By Lemma \ref{RRtriple}(iii), in this case $h^1(V,\calO_V(nH)) = 2$ for all $n\geq 1$. Using this, the remainder of the proof follows by the same argument used to prove \cite[Theorem 10.4]{bgd4a}.
\end{proof}

This completes the analysis of the elliptic ruled components.

\section{Type II Fibres} \label{type2sect}

In order to prove Theorem \ref{2polclass}, we plan to use \cite[Theorem 10]{bgd4a} and Theorem \ref{actriple} to study the induced linear systems on the components of the central fibre $X_0$. However, in order to apply these results we first need to ensure that certain assumptions are satisfied. In this section and the next we will show that we can achieve this by birationally modifying $X$ in a way that does not affect the form of its relative log canonical model.

We begin this section by setting up some notation. Let $\pi\colon X \to \Delta$ be a degeneration of K3 surfaces of degree two with polarisation divisor $H$ that satisfies Assumption \ref{mainass}, and suppose that $X_0 = \pi^{-1}(0)$ is a fibre of Type II. Write $X_0 = V_1 \cup \cdots \cup V_r$, with $V_2,\ldots,V_{r-1}$ elliptic ruled and $V_1$ and $V_r$ rational. Let $D_{i-1,i}$ denote the elliptic double curve ${V_{i-1} \cap V_i}$.

Then we have:

\begin{theorem} \label{type2pol} Let $\pi\colon X \to \Delta$ be a degeneration of K3 surfaces of degree two with polarisation divisor $H$ that satisfies Assumption \ref{mainass}, and assume that $X_0 = \pi^{-1}(0)$ is a fibre of Type II. Then with notation as above, no component of the double locus on $V_i$ is fixed in the linear system $|H_i|$.
\end{theorem}
\begin{proof} We proceed by contradiction. Let $H_i$ be the divisor on $V_i$ obtained by intersecting with $H$ and suppose that $D_{ij}$ is a double curve on $V_i$ (so $j \in \{i-1, i+1\}$) that is fixed in $|H_i|$. We aim to show that this implies that $(D_{ij}|_{V_i})^2 = 0$ and $H_i.D_{ij} = 0$. For then, in the notation of Section \ref{compdenfib}, $D_j|_{V_{ij}}$ is a nonsingular elliptic $0$-curve with self-intersection number $0$, which cannot exist by \cite[Lemma 2.2]{bgd4}.

So, in order to prove Theorem \ref{type2pol}, we need to show that $(D_{ij}|_{V_i})^2 = 0$ and $H_i.D_{ij} = 0$. The first of these will follow from the triple point formula \cite[Corollary 2.4.2]{odas} if we can show that $D_{ij}$ has non-positive self-intersection on both of the components in which it lies. The first step to proving this is to show that $D_{ij}$ is fixed in both of the components in which it lies.

\begin{remark} We note that the next few results are proved in considerably more generality than we need in order to prove Theorem \ref{type2pol}. However, the same results will also be used when we come to analyse the Type III fibres in the next section, and we will need the greater generality there. \end{remark}

\begin{lemma} \label{fixD} Let $V_i$ and $V_j$ be two distinct surfaces meeting along a double curve $D_{ij}$. Let $\calL$ be an invertible sheaf on $V = V_i \cup V_j$, such that there exist sections of $H^0(V, \calL)$ which are nonvanishing on $V_k$ for each $k$. Let $H_k$ denote an effective divisor on $V_k$ defined by a nonvanishing section of $H^0(V, \calL)$. Suppose that $D_{ij}$ is a fixed component of $\left| H_i \right|$ on $V_i$. Then $D_{ij}$ is also a fixed component of $\left| H_j \right|$ on $V_j$. 
\end{lemma}
\begin{proof} Assume first that $i \neq j$. For a contradiction, suppose that $D_{ij}$ is fixed in $\left| H_i \right|$ but not in $\left| H_j \right|$. Then there exists a section $s \in H^0(V, \calL)$ such that $s$ restricted to $V_j$ does not vanish on $D_{ij}$. But then $s$ restricted to $V_i$ defines a divisor linearly equivalent to $H_i$ that does not vanish on $D_{ij}$. However, this contradicts $D_{ij}$ being a fixed component of $\left| H_i \right|$. 
\end{proof}

Applying this lemma with $\calL$ equal to the restriction of $\calO_X(H)$ to $X_0$, and noting that this restriction defines the complete linear system $|H_i|$ on each component $V_i$ of $X_0$ by Lemma \ref{fibreres}, we see that $D_{ij}$ is fixed in both of the components in which it lies. The fact that it has non-positive self-intersection on both of these components will follow from another lemma:

\begin{lemma} \label{Dneg} Let $V_i$ be a normalised component of the central fibre $X_0$ in a degeneration of K3 surfaces $\pi\colon X\to \Delta$ of Type II or III, and let $D_i$ be the locus of double curves on $V_i$. Let $|H_i|$ be a linear system on ${V}_i$ which contains in its fixed locus an irreducible component $D_{ij}$ of ${D}_i$. Then $D_{ij}^2 \leq 0$, and this inequality is strict if $\pi\colon X\to \Delta$ is of Type III and $D_{ij}$ is smooth.
\end{lemma}
\begin{proof} Suppose for a contradiction that $D_{ij}^2 > 0$ (or $D_{ij}^2 \geq 0$ in the case where $\pi\colon X\to \Delta$ is of Type III and $D_{ij}$ is smooth). However, by \cite[Lemma 5]{bgd4a} (when $V_i$ is rational) or Lemma \ref{RRtriple} (when $V_i$ is elliptic ruled) $h^0({V}_i, \calO_{{V}_i}(D_{ij})) \geq 2$, so that $D_{ij}$ moves in a linear system of dimension $\geq 1$ and hence cannot be fixed.\end{proof}

Now the fact that $(D_{ij}|_{V_i})^2 = 0$ follows from Lemma \ref{fixD}, Lemma \ref{Dneg} and the triple point formula \cite[Corollary 2.4.2]{odas}. Therefore it just remains to show that $H_i.D_{ij} = 0$. This will follow from:

\begin{prop} \label{HD0} Let $V_k$ be the normalised components of the central fibre $X_0$ in a degeneration of K3 surfaces $\pi\colon X\to \Delta$ of Type II or III, and let $D_k$ denote the locus of double curves on $V_k$. Let $\calL$ be a nef line bundle on $X_0$ such that $\calL.\calL = 2$ and there exist nonvanishing sections in $H^0(V_k, \calL)$ for all $k$. Let $H_k$ denote a divisor defined on ${V}_k$ by such a section. Then if an irreducible component ${D}_{ij} = V_i \cap V_j$ of $D_i$ is in the fixed locus of $|H_i|$ for some $i$, it must satisfy $H_i.{D}_{ij} = 0$. \end{prop}

\begin{remark} We remark that this proposition is the only place in the proof of Theorem \ref{2polclass} that the degree two assumption (that $H_t^2 = 2$ for all $t \in \Delta^*$) is really essential; as noted in Remark \ref{ellrem}, the results of Subsection \ref{ellcomp} which also make this assumption should admit generalisations to higher degrees. Unfortunately, we expect that the conclusions of Proposition \ref{HD0} will be false in much higher degrees, although it seems feasible that the proof below may admit a generalisation to the degree four case. \end{remark}

\begin{proof}[Proof of Proposition \ref{HD0}] Before we begin with the proof, we make a remark about non-normal components. If $V_i$ is a non-normal component in $X_0$ it intersects itself along a smooth rational curve $D_{ii}$. When we normalise it we find that $D_{ii}$ has two preimages. As these preimages will usually be considered alongside double curves that lie in two different components, we will abuse notation and refer to them as $D_{ij}|_{V_i}$ and $D_{ij}|_{V_j}$, where it is understood that if $i$ and $j$ are equal these refer to the disjoint curves in the normalisation. Finally, we note that in this case if one of these curves is in the fixed locus of $|H_i|$ then the other must also be.

We now proceed with the proof. Suppose $D_{ij}$ is in the fixed locus of $|H_i|$. Then, using Lemma \ref{fixD} if $i\neq j$, we see that $D_{ij}$ is also in the fixed locus of $|H_j|$. So, by Lemma \ref{Dneg}, $(D_{ij}|_{V_i})^2 \leq 0$ and $(D_{ij}|_{V_j})^2 \leq 0$, and these inequalities are strict when $\pi\colon X \to \Delta$ is Type III and $D_{ij}$ is smooth. Putting this into the triple point formula \cite[Corollary 2.4.2]{odas}
\begin{equation*}(D_{ij}|_{V_i})^2 + (D_{ij}|_{V_j})^2 = \begin{cases} 0 & \mathrm{if\ Type\ II\ or\ Type\ III\ with\ } D_{ij}\ \mathrm{nodal} \\ -2 & \mathrm{if\ Type\ III\ with\ } D_{ij}\ \mathrm{smooth} \end{cases} \end{equation*} 
we get that $(D_{ij}|_{V_i})^2 = (D_{ij}|_{V_j})^2 = 0$ if $\pi\colon X \to \Delta$ is of Type II or Type III with $D_{ij}$ nodal, and $(D_{ij}|_{V_i})^2 = (D_{ij}|_{V_j})^2 = -1$ if $\pi\colon X \to \Delta$ is of Type III with $D_{ij}$ smooth.

Now write the linear system $\left| H_i \right|$ as
\begin{equation*} |H_i| = |H_i|_m + |H_i|_f, \end{equation*}
where $|H_i|_f$ is the fixed part of $|H_i|$ and $|H_i|_m$ has no fixed components. Note that $|H_i|_m$ and $|H_i|_f$ are effective or trivial and $|H_i|_m$ is nef. 

Suppose that we are on a component $V_i$ with $H_i^2 = 0$. Then
\begin{equation*} 0 = H_i^2 = |H_i|_m^2 + |H_i|_m.|H_i|_f + H_i.|H_i|_f. \end{equation*}
As $H_i$ and $|H_i|_m$ are nef, all terms on the right hand side of this equation are zero. Furthermore, by effectiveness of $|H_i|_f$, we have that $H_i.F = 0$ for all irreducible components $F$ of $|H_i|_f$. Hence if $D_{ij}$ is fixed in $|H_i|$, then $H_i.D_{ij}= 0$.

So we are left with the case where $D_{ij}$ is the intersection of two components $V_i$ and $V_j$ with $H_i^2 >0$ and $H_j^2>0$. As $\calL.\calL = 2$ on $X_0$, this can only occur if $H_i^2 = H_j^2 = 1$, or if $i = j$ and $H_i^2 = 2$ (in which case $V_i$ is the normalisation of a surface that intersects itself). In these cases we can explicitly analyse the form of $|H_i|$ on $V_i$.

Let $|H_i|$ be a linear system on $V_i$, with $H_i^2 = 1$ or $2$ and $D_{ij}$ be a double curve that is a fixed component of $|H_i|$. As above, write
\begin{equation} H_i^2 = |H_i|_m^2 + |H_i|_m.|H_i|_f + H_i.|H_i|_f. \label{fixedH} \end{equation}
If $H_i.|H_i|_f = 0$, we are done as above. So assume $H_i.|H_i|_f > 0$. We will show that this implies that $|H_i|_f^2 > 0$ and $|H_i|_m.|H_i|_f = 0$, then use these expressions to derive a contradiction. 

If $H_i^2 = 1$, then $H_i.|H_i|_f = 1$ and so
\begin{equation*}1 = H_i.|H_i|_f = |H_i|_m.|H_i|_f + |H_i|_f^2.\end{equation*}
As $|H_i|_m$ is nef, by equation \eqref{fixedH} necessarily $|H_i|_m.|H_i|_f = 0$. So $|H_i|_f^2 = 1$. 

If $H_i^2 = 2$ then necessarily $i=j$, so $D_{ij}|_{V_i}$ and $D_{ij}|_{V_j}$ both lie in $V_i$ and so are both in $|H_i|_f$. In this case, write 
\begin{equation*}|H_i|_f = a_0 (D_{ij}|_{V_i}) + a_1 (D_{ij}|_{V_i}) + \sum_{k =2}^m a_k F_k,\end{equation*}
for irreducible curves $F_k$ and integers $a_k$ with $a_k > 0$. Then
\begin{equation*}H_i.|H_i|_f = a_0 H_i.(D_{ij}|_{V_i}) + a_1 H_i.(D_{ij}|_{V_j}) + \sum_{k =2}^m a_k H_i.F_k.\end{equation*}
As $H_i$ is nef, all of the terms in the right hand side of this equation are non-negative. Furthermore, we may assume $H_i.(D_{ij}|_{V_i})>0$, otherwise we are done. But $H_i.(D_{ij}|_{V_i}) =  H_i.(D_{ij}|_{V_j})$, so $H_i.|H_i|_f \geq 2$. Therefore, by equation \eqref{fixedH} and the fact that $|H_i|_m$ is nef, $H_i.|H_i|_f = 2$ and $|H_i|_m.|H_i|_f = 0$. Then as 
\begin{equation*}2 = H_i.|H_i|_f = |H_i|_m.|H_i|_f + |H_i|_f^2,\end{equation*}
we must have $|H_i|_f^2 = 2$.

In either case $|H_i|_f^2 > 0$ and $|H_i|_m.|H_i|_f =0$. Note further that, for any component $F$ of $|H_i|_f$, we must have $|H_i|_m.F = 0$ as $|H_i|_f$ is effective and $|H_i|_m$ is nef. So, for all components $F$ of $|H_i|_f$, 
\begin{equation*}|H_i|_f.F = |H_i|_m.F + |H_i|_f.F = H_i.F \geq 0\end{equation*}
as $H_i$ is nef. Therefore as $|H_i|_f$ is effective, it must be nef.

Now consider the case where $V_i$ is elliptic ruled. Then $|H_i|_f^2 = 1$, as no component may intersect itself in a degeneration of Type II. Let $\hat{V}_i$ be the surface obtained from $V_i$ by contracting all rational $(-1)$-curves $C$ on $V_i$ with $|H_i|_f.C = 0$, and let $|\hat{H}_i|_f$ and $\hat{D}_{ij}$ denote the images on $\hat{V}_i$ of $|H_i|_f$ and $D_{ij}$ respectively. Then $\hat{V}_i$ is $|\hat{H}_i|_f$-minimal with $|\hat{H}_i|_f.\hat{D}_{ij} = |H_i|_f.D_{ij}$ and $(\hat{D}_{ij})^2 \geq (D_{ij}|_{V_i})^2$ for each choice of $j \in \{i-1,i+1\}$.

Assume first that $(\hat{D}_{i,i-1})^2$ and $(\hat{D}_{i,i+1})^2$ are both zero. Then Lemma \ref{HD'0} shows that $\hat{V}_i$ is minimally ruled and $|\hat{H}_i|_f.\hat{D}_{i,i-1} = |\hat{H}_i|_f.\hat{D}_{i,i+1} > 0$. So we also have $|H_i|_f.D_{i,i-1} = |H_i|_f.D_{i,i+1} > 0$ and Lemma \ref{RRtriple}(i) applied to $|H|_f$ gives $h^0(V_i, \calO_{V_i}(|H_i|_f)) \geq 2$, contradicting $|H_i|_f$ being fixed.

Next, consider the case where $(\hat{D}_{i,i-1})^2$ and $(\hat{D}_{i,i+1})^2$ are not both zero. Then, using Lemma \ref{HD'0}, we may assume that $|H_i|_f.D_{ik} = |\hat{H}_i|_f.\hat{D}_{ik} = 0$ and that $(D_{ik}|_{V_i})^2 \leq (\hat{D}_{ik})^2 < 0$ for some $k \in \{i-1,i+1\}$. Thus, as any fixed double curve in a Type II degeneration has $(D_{ik}|_{V_i})^2 =0$, we see that $D_{ik}$ cannot be in $|H_i|_f$. But then, using Lemma \ref{RRtriple} and the fact that $|H_i|_f.D_{ik}  = 0$, we see that $h^0(V_i, \calO_{V_i}(|H_i|_f)) \geq 2$, contradicting $|H_i|_f$ being fixed.

Finally, consider the case where $V_i$ is rational. Then $|H_i|_f$ is effective, $|H_i|_f^2 > 0$ and $|H_i|_f.D_i \geq 0$ as $|H_i|_f$ is nef. So \cite[Lemma 5]{bgd4a} gives that $h^0(V_i, \calO_{V_i}(|H_i|_f)) \geq 2$, contradicting $|H_i|_f$ being fixed. \end{proof}

Applying this proposition with $\calL$ equal to the restriction of $\calO_X(H)$ to $X_0$ gives $H_i.D_{ij} = 0$, as required. This completes the proof of Theorem \ref{type2pol}. \end{proof}

\section{Type III Fibres}

We next consider the case where $\pi\colon X \to \Delta$ is of Type III. In the same way as Section \ref{type2sect}, we would like to show that we may birationally modify $X$ so that the assumptions of \cite[Theorem 10]{bgd4a} and Theorem \ref{actriple} hold.  

As before, we begin by setting up some notation. Let $\pi\colon X \to \Delta$ be a degeneration of K3 surfaces of degree two with polarisation divisor $H$ that satisfies Assumption \ref{mainass}, and assume that $X_0 = \pi^{-1}(0)$ is a fibre of Type III. Write $X_0 = V_1 \cup \cdots \cup V_r$ for rational surfaces $V_i$. Note that we assume here that the surfaces $V_i$ have been normalised. Let $D_{ij}$ denote the rational double curve $V_{i} \cap V_j$ and let $D_i = \bigcup_{j} D_{ij}$ denote the double locus on $V_i$. We call a double curve $D_{ij}$ a ($*$)-\emph{curve} if it is nonsingular and has self-intersection $(-1)$ on both components in which it lies. Finally, let $H_i$ denote the effective (or zero) divisor on $V_i$ defined by intersecting with $H$.

\begin{remark} Before we embark on the results of this section, we make a remark on non-normal components. Suppose that $\overline{V}_i$ is a non-normal component in the Type III fibre $X_0$. Let $\nu\colon V_i \to \overline{V}_i$ denote the normalisation. Then the non-normal locus in $\overline{V}_i$ is a smooth rational curve $\overline{D}_{ii}$. The preimage $D_{ii} := \nu^{-1}(\overline{D}_{ii})$ consists of two disjoint rational curves in $V_i$. As these curves will often be considered alongside double curves that lie in two different components, we will frequently abuse notation and refer to them as $D_{ij}|_{V_i}$ and $D_{ij}|_{V_j}$, where it is understood that if $i$ and $j$ are equal these refer to the disjoint curves in the normalisation. Finally, when we refer to the restriction of a line bundle or divisor to a component $V_i$, we understand this to mean the pull-back of that intersection under the normalisation map $\nu$. \end{remark}

Now we can state the main result of this section.

\begin{theorem} \label{type3pol} Let $\pi\colon X \to \Delta$ be a degeneration of K3 surfaces of degree two with polarisation divisor $H$ that satisfies Assumption \ref{mainass}, and assume that $X_0 = \pi^{-1}(0)$ is a fibre of Type III. Then with notation as above, if a component $D_{ij}$ of the double locus $D_i$ on $V_i$ is a fixed component of $\left| H_i \right|$, then $D_{ij}$ is a $(*)$-curve with $H_i.D_{ij} = 0$. \end{theorem}
\begin{proof} Suppose that $D_{ij} = V_i \cap V_j$ is fixed in the linear system $|H_i|$. Applying Proposition \ref{HD0} with $\calL$ equal to the restriction of $\calO_X(H)$ to $X_0$, and noting that this restriction defines the complete linear system $|H_i|$ on each component $V_i$ of $X_0$ by Lemma \ref{fibreres}, we get that $H_i.D_{ij} = 0$, proving half of the theorem. 

Next, using Lemma \ref{fixD} we see that $D_{ij}$ must also be fixed in $|H_j|$ (or, if $i=j$, we note that both components of $D_{ii}$ must be fixed in $|H_i|$). This places a strong restriction on $D_{ij}$: by Lemma \ref{Dneg} it must have negative self-intersection in both $V_i$ and $V_j$ (or non-positive if $D_{ij}$ is nodal). In this setting, the triple point formula \cite[Corollary 2.4.2]{odas} gives
\begin{equation*}(D_{ij}|_{V_i})^2 + (D_{ij}|_{V_j})^2 = \begin{cases} 0 & \mathrm{if} \ D_{ij} \ \mathrm{is}\ \mathrm{a}\ \mathrm{nodal}\ \mathrm{curve}\ \mathrm{on}\ V_i\ \mathrm{or}\ V_j \\ -2 & \mathrm{otherwise} \end{cases} \end{equation*}
If $D_{ij}$ is nonsingular, this implies that it must be a $(*)$-curve, as required. So it only remains to show that $D_{ij}$ cannot be a rational nodal curve. However, if $D_{ij}$ is a rational nodal curve on $V_i$ then the argument above shows that it is a double curve that is both a Kodaira elliptic curve and a $0$-curve on $V_i$. Such curves cannot exist by \cite[Lemma 2.2]{bgd4}. This completes the proof of Theorem \ref{type3pol}. \end{proof}

\section{Proof of Theorem \ref{2polclass}}

We now embark upon the proof of Theorem \ref{2polclass}.  

As in the statement of Theorem \ref{2polclass}, let $\pi\colon X \to \Delta$ be a degeneration of K3 surfaces of degree two with polarisation divisor $H$ that satisfies Assumption \ref{mainass}. Then by Lemma \ref{fibreres}, we may safely restrict our attention to the pair $(X_0,H_0)$, as the natural morphism $\phi_0\colon X_0 \to (X_0)^c$ to the log canonical model of $(X_0,H_0)$ agrees with the restriction of $\phi$ to $X_0$. 

Now, as in the previous sections, let $V_i$ denote a normalised irreducible component of $X_0$ and let $H_i$ denote the effective divisor induced on $V_i$ by $H$. The linear systems $|nH_i|$ for $n \geq 0$ define a morphism
\begin{equation*}\phi_{V_i}\colon V_i \longrightarrow \mathrm{Proj}\bigoplus_{n \geq 0}H^0(V_i, \calO_{V_i}(nH_i)).\end{equation*}
The restriction of $\phi_0$ to $V_i$ factors through this morphism. So we can gather information about $\phi_0$ by studying the maps $\phi_{V_i}$.

We split our study of the maps $\phi_{V_i}$ into three cases, depending upon whether $V_i$ is a $0$-, $1$- or $2$-surface.

\begin{lemma} \label{type3cont} If $V_i$ is a $1$-surface \textup{(}resp. $2$-surface\textup{)}, then $V_i$ is contracted to a curve \textup{(}resp. point\textup{)} by $\phi_{V_i}$.  \end{lemma}
\begin{proof} If $V_i$ is a $1$-surface, Lemma \ref{012facts} gives that the pencil of $0$-curves on $V_i$ form a ruling. $\phi_{V_i}$ contracts all $0$-curves on $V_i$, so $\phi_{V_i}$ contracts $V_i$ onto a section of the ruling.

If $V_i$ is a $2$-surface, then $H_i$ is numerically trivial. So $H_i.E = 0$ for any effective divisor $E$ on $V_i$, and $\phi_{V_i}$ contracts $V_i$ to a point. \end{proof}

Given this, we see that the form of $(X_0)^c$ depends mostly upon what happens to the $0$-surfaces. The crux of Theorem \ref{2polclass} is the following result:

\begin{theorem} \label{0surfprop} After performing a birational modification of $X_0$ that does not affect the form of $\phi_0(X_0)$, for any $0$-surface $V_i$ we may assume that the linear system $|nH_i|$ has no fixed components or base locus for $n \geq 2$ and that the morphism $\phi_{nH_i}$ to projective space defined by this linear system is birational onto its image for $n \geq 3$. 
\end{theorem}

\begin{proof} First suppose that $X_0$ is a fibre of Type I. Then $X_0$ is a smooth K3 surface and $H_0$ is nef and big on $X_0$. In this case Mayer \cite[Proposition 8]{fk3s} shows that $|H_0|$ is either base point free or has the form $|H_0| = |2E| + F$, where $E$ is a nonsingular elliptic curve with $E^2 = 0$ and $F$ is a fixed rational curve with $F^2 = -2$ and $E.F = 1$. In either case, the proof of \cite[Corollary 5]{fk3s} shows that $|2H_0|$ has no fixed components and $\phi_{3H_0}$ is birational onto its image.

In the cases where $X_0$ is a fibre of Type II or III, the theorem will follow from \cite[Theorem 10]{bgd4a} and Theorem \ref{actriple} once we have used Theorems \ref{type2pol} and \ref{type3pol} to prove that we may assume that the conditions of these results hold.

So assume that $X_0$ is a fibre of Type II. Suppose first that $V_i$ is a rational $0$-surface in $X_0$. By Theorem \ref{type2pol} we see that the linear system $|H_i|$ does not contain any component of $D_i$ in its fixed locus, so we may apply \cite[Theorem 10.4]{bgd4a} to prove Theorem \ref{0surfprop}.

Next suppose that $V_i$ is an elliptic ruled $0$-surface. By Theorem \ref{type2pol} we see that the linear system $|H_i|$ does not contain any component of $D_i$ in its fixed locus. Define a birational modification of $X_0$ by contacting any rational $(-1)$-curves $C$ on $V_i$ with $H_i.C = 0$. This does not affect the nonsingularity of $V_i$ or the form of $\phi_0(X_0)$. The resulting surface is $H_i$-minimal and does not contain any component of $D_i$ in its fixed locus. So, by Lemma \ref{HD'0}, we may assume that $H_i.D_{ij} = 0$ for one of the double curves $D_{ij}$ on $V_i$. Therefore, the assumptions of Theorem \ref{actriple} are satisfied and we may apply this result to prove Theorem \ref{0surfprop}. 

Finally, assume that $X_0$ is a fibre of Type III. Then every $0$-surface $V_i$ is rational. In order to apply \cite[Theorem 10.4]{bgd4a} in this case, we need to show that we may assume that no component of the double curve $D_i$ is fixed in $|H_i|$. In order to show this, we will use Theorem \ref{type3pol} to prove that we may contract any such $D_i$, and that this contraction does not affect the form of $\phi_0(X_0)$ or the nonsingularity of the $0$-surfaces on $X_0$ (although $X_0$ itself may become non-semistable).

So suppose that $V_i$ is a $0$-surface in $X_0$ and let $D_{ij}$ be a component of the double curve in $V_i$ that is fixed in $|H_i|$. Then, by Theorem \ref{type3pol}, $D_{ij}$ is a $(*)$-curve and $H_i.D_{ij} = 0$. Perform a Type II modification to move $D_{ij}$ to a neighbouring component. This has the effect of contracting $D_{ij}$ on $V_i$, but does not affect the assumptions of Theorem \ref{type3pol}. So we may repeat the argument for all components of the double curve on $V_i$ that are fixed in $|H_i|$. The result is that all such $D_{ij}$ are contracted on $V_i$ without affecting its nonsingularity.

Now define a birational modification of $X_0$ by contacting any double curve that is fixed in the linear system $|H_i|$ on a $0$-surface $V_i$. By the argument above, this does not affect the nonsingularity of the $0$-surfaces (although it may lead to $X_0$ becoming non-semistable). Furthermore, as every such curve has $H_0.D_{ij} = 0$, it does not affect the form of $\phi_0(X_0)$ either. Finally, after performing this birational modification we see that no $0$-surface $V_i$ contains a double curve in the fixed locus of $|H_i|$, so we may use \cite[Theorem 10.4]{bgd4a} to complete the proof of Theorem \ref{0surfprop} in this case.\end{proof}

\begin{proof}[Proof of Theorem \ref{2polclass}] First consider the case where $\pi\colon X \to \Delta$ is a degeneration of Type I. Then $X_0$ is a smooth K3 surface. We split into two cases, depending upon whether $|H_0|$ is base point free or not. If $|H_0|$ is base point free, Theorem \ref{0surfprop} and a Riemann-Roch calculation shows that $\phi_0(X_0)$ is hyperelliptic. If $|H_0|$ has base points, a similar calculation shows that $\phi_0(X_0)$ is unigonal.

Next assume that $\pi\colon X \to \Delta$ is a degeneration of Type II or III. By Lemma \ref{type3cont} and Theorem \ref{0surfprop} we see that $\phi_0$ contracts all $1$- and $2$-surfaces to curves and points, and maps $0$-surfaces birationally onto their images. Furthermore, as $H_i^2 > 0$ on a $0$-surface $V_i$ (by Lemma \ref{Hpos}) and $H_0^2 = 2$, there are at most two $0$-surfaces in $X_0$. Thus, we may assume that the contraction of the $1$- and $2$-surfaces has taken place and restrict to the case where $X_0$ is a union of at most two irreducible components. We split into subcases depending upon the number and nature of these components:
\begin{list}{}{}
\item[(A)] There are two rational components, $V_1$ and $V_2$.
\item[(B)] There is only one component, $V_1$, that is rational.
\item[(C)] There are two components, $V_1$ and $V_2$, where $V_1$ is rational and $V_2$ is elliptic ruled.
\item[(D)] There are two elliptic ruled components, $V_1$ and $V_2$.
\item[(E)] There is only one component, $V_1$, that is elliptic ruled.
\end{list}
Given this, Theorem \ref{2polclass} is proved by using Theorem \ref{0surfprop} and the Riemann-Roch theorem to analyse each of the above cases and show that they have the required type. The full analysis is rather long and will not be reproduced here; we refer the interested reader to \cite[Chapter 3]{mythesis}. Instead, we will briefly analyse each of the above cases in turn and state how they correspond to the cases listed in Tables \ref{eqtable1} and \ref{eqtable2}.

Let $D$ denote the double locus on $X_0$. Assume first that we are in case (A). By the genus formula, there are three distinct possibilities, (A1), (A2) and (A3), corresponding to when $(H_1.D,H_2.D) = (1,1)$, $(3,3)$ or $(1,3)$ respectively. Given this, an easy Riemann-Roch calculation shows that case (A1) corresponds to cases II.4 (where there is no $\tilde{E}_8$ singularity present) or III.4 (where the locus $\{f_6 = l_i = 0\}$ is either reduced for both $i$ or non-reduced for both $i$) in Table \ref{eqtable2}, depending upon whether $X_0$ is of Type II or III respectively. Similarly, case (A2) corresponds to II.3 or III.3 in Table \ref{eqtable1}. Case (A3) can only happen if the degeneration is of Type III and $V_2$ intersects itself along a component of $D$. In this situation we again find ourselves in case III.4 of Table \ref{eqtable2}, but this time the locus $\{f_6 = l_i = 0\}$ is allowed to have a component of multiplicity $2$ for one choice of $i \in \{1,2\}$. We remark here that, in the Type II case, the calculations in cases (A1) and (A2) are completely analogous to Friedman's calculations in cases (5.2.2) and (5.2.1) in \cite{npgttk3s} respectively.

Before we can analyse case (B), we need a lemma that will allow us to treat the case where $\pi\colon X \to \Delta$ is a degeneration of Type III and, prior to contracting the other components, the polarisation $H_1$ on $V_1$ satisfied $H_1.D = 0$.

\begin{lemma} \label{1cyclelem} Let $\pi\colon X \to \Delta$ be a degeneration of K3 surfaces of degree two with polarisation divisor $H$ that satisfies Assumption \ref{mainass}, and assume that $X_0 = \pi^{-1}(0)$ is a fibre of Type III. Then if $V_i$ is a $0$-surface in $X_0$ with $H_i.D_{ij} = 0$ for all double curves $D_{ij}$ in $V_i$, 
\begin{list}{\textup{(\alph{list2})}}{\usecounter{list2}}
\item $V_i$ is the only $0$-surface in $X_0$, and
\item all of the other components of $X_0$ are $2$-surfaces.
\end{list}\end{lemma}
\begin{proof} We begin with (a). Note that $H_i^2 > 0$ by Lemma \ref{Hpos}, so $H_i^2$ must equal 1 or 2. Then, by the Riemann-Roch Theorem, $( H_i^2 + \sum_j H_i.D_{ij}) \in 2\Z$. If $H_i^2 = 1$ then $H_i.D_{ij} > 0$ for some $j$, as $H_i$ is nef, and we are done. So we may assume $H_i^2 = 2$. By Lemma \ref{Hpos} again, $V_i$ can be the only $0$-surface in $X_0$.

Now consider (b). We want to show that $X_0$ cannot contain any $1$-surfaces. Suppose that this is not the case, so that there exists a $1$-surface $V_{j_1}$ in $X_0$. By Proposition \ref{012facts} we see that $V_{j_1}$ is ruled by the pencil of $0$-curves on it. Let $F_{j_1}$ be a general fibre in this ruling that is not a component of the double curve $D_{j_1}$. Then, by \cite[Lemma 3]{bgd4a}, $F_{j_1}.D_{j_1} = 2$. As $D_{j_1}$ is effective, this implies that $D_{j_1}$ contains either two sections or one bisection of the ruling for $V_{j_1}$. Suppose that $D_{j_1j_2} = V_{j_1} \cap V_{j_2}$ is one such section or bisection.

Now consider the divisor $H_{j_1}$ on $V_{j_1}$ defined by intersecting with $H$. By Lemma \ref{Hpos} we have $H_{j_1}^2 = 0$ so, as $H_{j_1}$ is effective and nef, $H_{j_1}$ must be a sum of $0$-curves. Thus $H_{j_1}$ is a sum of fibres of the ruling for $V_{j_1}$ and $H_{j_1}.D_{j_1j_2} > 0$ as $D_{j_1j_2}$ is a section or bisection of this ruling. 

Next consider the component $V_{j_2}$. As $H_{j_2}.D_{j_1j_2}= H_{j_1}.D_{j_1j_2} > 0$, this component cannot be a $2$-surface. Moreover it cannot be a $0$-surface, as by assumption $V_i$ is the only $0$-surface and $H_i.D_{ik} = 0$ for all $k$. Therefore, $V_{j_2}$ must be another $1$-surface. 

Repeating the argument above, we see that $H_{j_2}$ must be a sum of fibres of a ruling for $V_{j_2}$ and $D_{j_1j_2}$ must be a section or bisection of that ruling. Furthermore, if $D_{j_1j_2}$ is a section then there is another double curve $D_{j_2j_3}$ on $V_{j_2}$ that is a section of the same ruling, so we may repeat the process to find a $1$-surface $V_{j_3}$ meeting $V_{j_2}$ along $D_{j_2j_3}$. If $D_{j_1j_2}$ is a bisection then no other double curves on $V_{j_2}$ are sections or bisections of the ruling, so the process terminates.

Repeating this argument as many times as possible and relabeling components if necessary, we obtain either:
\begin{itemize}
\item a cycle $V_{j_1},\ldots,V_{j_n}$ of $1$-surfaces meeting along sections of a given ruling for each; or
\item a chain $V_{j_1},\ldots,V_{j_n}$ of $1$-surfaces such that $D_{j_kj_{k+1}} = V_{j_k} \cap V_{j_{k+1}}$ is a section of a given ruling on $V_{j_k}$ for $k \in \{2,\ldots,n-1\}$ and a bisection of a given ruling on $V_1$ and $V_n$.
\end{itemize}
Moreover $H_{j_k}.F_{j_k} = 0$ for any fibre $F_{j_k}$ of the given ruling on $V_{j_k}$ and any $1 \leq k \leq n$. However, such configurations of $1$-surfaces are excluded by \cite[Lemma 2.2]{bgd4}. This is a contradiction, so $X_0$ cannot contain any $1$-surfaces. 
\end{proof}

Now we are ready to consider case (B). In this case there is only one component, but this component may be non-normal with self-intersection along $D$. By the genus formula, there are three distinct possibilities, (B1), (B2) or (B3), corresponding to the cases where $H_1.D = 2$, $H_1.D = 4$ or $D = \emptyset$ respectively. Another easy calculation shows that case (B1) corresponds to II.1 (where there is no $\tilde{E}_7$ singularity present) or III.1 and case (B2) corresponds to II.2 or III.2 in Table \ref{eqtable1}. In subcase (B3), \cite[Proposition 2.5]{bgd4} and Lemma \ref{1cyclelem} (in the Type II and III cases respectively) show that, prior to contracting the other components, $V_1$ was the only $0$-surface in $X_0$ and all other components were $2$-surfaces. These $2$-surfaces were contracted to an elliptic or cusp singularity. By using Theorem \ref{0surfprop} and carefully analysing the curves that were contracted to the singularity on $V_1$, we see that subcase (B3) corresponds to II.0h (where there is exactly one $\tilde{E}_7$ or $\tilde{E}_8$ singularity present) or III.0h in Table \ref{eqtable1} if $|H_1|$ is base point free and II.0u (where there is exactly one $\tilde{E}_7$ or $\tilde{E}_8$ singularity present) or III.0u if $|H_1|$ has base points. Once again we note here that, in the Type II case, the calculations in cases (B1) and (B2) are completely analogous to Friedman's calculations in cases (5.2.3) and (5.2.4) in \cite{npgttk3s} respectively.

Cases (C), (D) and (E) can only arise if our degeneration is of Type II. In order to analyse them, we require one final lemma.

\begin{lemma} \label{type20surf} Let $\pi\colon X \to \Delta$ be a degeneration of K3 surfaces of degree two with polarisation divisor $H$ that satisfies Assumption \ref{mainass}, and assume that $X_0 = \pi^{-1}(0)$ is a fibre of Type II. If $V_i$ is an elliptic ruled $0$-surface in $X_0$, then $V_i$ must intersect a $2$-surface.
\end{lemma}
\begin{proof} Let $\hat{V}_i$ be the surface obtained from $V_i$ by contracting all rational $(-1)$-curves $C$ with $H_i.C = 0$, and let $\hat{H}_i$ and $\hat{D}_{ij}$ denote the images on $\hat{V}_i$ of $H_i$ and $D_{ij}$ respectively. Then $\hat{V}_i$ is $\hat{H}_i$-minimal with $\hat{H}_i^2 = H_i^2$, $\hat{H}_i.\hat{D}_{ij} = H_i.D_{ij}$ and $(\hat{D}_{ij})^2 \geq (D_{ij}|_{V_i})^2$.

Assume first that $(\hat{D}_{i,i-1})^2$ and $(\hat{D}_{i,i+1})^2$ are not both zero. Then, by Lemma \ref{HD'0}, we may assume that $H_i.D_{ik} = \hat{H}_i.\hat{D}_{ik} = 0$ and $(D_{ik}|_{V_i})^2 \leq (\hat{D}_{ik})^2 < 0$ for some $k \in \{i-1,i+1\}$. Since $D_{ik} = V_i \cap V_{k}$, by the triple point formula \cite[Corollary 2.4.2]{odas}, $(D_{ik}|_{V_{k}})^2 > 0$. But then $D_{ik}$ is an effective divisor on $V_{k}$ with strictly positive self-intersection and $H_k.D_{ik} = 0$ so, by Proposition \ref{012facts}, $V_{k}$ is a $2$-surface.

Next, suppose that $(\hat{D}_{i,i-1})^2$ and $(\hat{D}_{i,i+1})^2$ are both zero. Then Lemma \ref{HD'0} shows that $\hat{V}_i$ is minimally ruled and $\hat{H}_i.\hat{D}_{i,i-1} = \hat{H}_i.\hat{D}_{i,i+1} > 0$. Analysing the class of $\hat{H}_i$ in $H^2(\hat{V}_i,\Z)$ using \cite[Proposition III.18]{cas}, we find that we must have $[\hat{H}_i] = s+f$, where $s$ is the class of a section and $f$ is the class of a fibre. This gives $\hat{H}_i^2 = 2$ and $\hat{H}_i.\hat{D}_{ij} = 1$ for any $j \in \{i-1,i+1\}$. So $H_i^2 = 2$ and $H_i.D_{ij} = 1$ for both values of $j$ also.

Therefore, by Lemma \ref{Hpos}, $V_i$ is the only $0$-surface. As $H_i.D_{ij} = 1$ for any choice of $j \in \{i-1,i+1\}$, $H_i$ cannot be numerically trivial on either of the components intersecting $V_i$, so both of these components must be $1$-surfaces. Therefore, by \cite[Proposition 2.5]{bgd4}, $X_0$ cannot contain any $2$-surfaces, so all surfaces except $V_i$ must be $1$-surfaces. Furthermore, using Proposition \ref{012facts} it is easily seen that if $V_j$ is any elliptic ruled $1$-surface, then $H_j$ must be a sum of fibres of the ruling on $V_j$.

Now consider $V_1$. By the argument above, $V_1$ is a rational $1$-surface with $H_1^2~=~0$. As all components between $V_1$ and $V_i$ are $1$-surfaces which are polarised by sums of fibres, we obtain that \mbox{$H_1.D_{1,2} = H_i.D_{i-1,i} = 1$}, where $D_{1,2}$ denotes the unique elliptic double curve on $V_1$.  Finally, $(V_1, D_{1,2})$ is an anticanonical pair, so by definition $K_{V_1}~=~-D_{1,2}$ and thus \mbox{$H_1.(H_1 + K_{V_1}) = -1$}. But the genus formula for effective divisors implies that this number must be even. This is a contradiction, so this case cannot occur and $V_i$ must intersect a $2$-surface.\end{proof}

Now consider case (C). Lemma \ref{type20surf} and \cite[Proposition 2.5]{bgd4} show that, prior to contracting the other components, $V_2$ met a chain of $2$-surfaces. These $2$-surfaces were contracted to give an elliptic singularity, which analysis using Lemma \ref{HD'0} shows to have type $\tilde{E}_8$. This case corresponds to II.4 in Table \ref{eqtable2}, where there is one $\tilde{E}_8$ singularity present.

Next consider case (D). Lemma \ref{type20surf} and \cite[Proposition 2.5]{bgd4} show that, prior to contracting the other components, $V_1$ and $V_2$ both met chains of $2$-surfaces. These $2$-surfaces were contracted to give one elliptic singularity in each of $V_1$ and $V_2$, which analysis using Lemma \ref{HD'0} shows to have type $\tilde{E}_8$. Thus, this case corresponds to II.4 in Table \ref{eqtable2}, where there are two $\tilde{E}_8$ singularities present. 

Finally, consider case (E). In this case Lemma \ref{HD'0} shows that, prior to contracting the other components, $V_1$ contained two double curves $D'$ and $D''$ with $H_1.D'=0$ and $H_1.D'' \geq 0$. We obtain two possibilities, (E1) and (E2), corresponding to $H_1.D'' = 2$ or $0$ respectively. By \cite[Proposition 2.5]{bgd4}, in case (E1) we see that, prior to contracting the other components, $V_1$ met a chain of $2$-surfaces along $D'$. This chain was contracted to give an elliptic singularity of type $\tilde{E}_7$. Furthermore, by \cite[Proposition 2.5]{bgd4} again, we see that $V_1$ met a chain of $1$-surfaces along $D''$ which were contracted onto the image of $D''$. This case corresponds to II.1 in Table \ref{eqtable1}, where there is an $\tilde{E}_7$ singularity present. Finally, in case (E2), \cite[Proposition 2.5]{bgd4} gives that, prior to contracting the other components, $V_1$ met two chains of $2$-surfaces along $D'$ and $D''$. These were contracted to give two elliptic singularities of type $\tilde{E}_8$. This gives case II.0h in Table \ref{eqtable2} (with two $\tilde{E}_8$ singularities) if $|H_1|$ is base point free, or case II.0u in Table \ref{eqtable2} (with two $\tilde{E}_8$ singularities) if $|H_1|$ has base points.\end{proof}

\bibliography{ThesisBooks}
\bibliographystyle{amsplain}

\end{document}